\documentclass[10pt]{amsart}

\usepackage{tikz}
\usepackage{amssymb,amsmath,latexsym,graphicx}
\usepackage{charter,eucal}
\usepackage{amssymb}
\usepackage[all,cmtip]{xy}
\usepackage{rotating,multicol}
\usetikzlibrary{decorations.markings}

\newtheorem{theorem}{Theorem }[section]

\newtheorem{lemma}[theorem]{Lemma}
\newtheorem{observation}[theorem]{Observation}

\newtheorem{remark}[theorem]{Remark}
\newtheorem{corollary}[theorem]{Corollary}
\newtheorem{proposition}[theorem]{Proposition}
\newtheorem{principle}[theorem]{\textsc{Principle}}

\newcommand{\bt}{\begin{theorem}}
\newcommand{\et}{\end{theorem}}
\newcommand{\bmt}{\begin{maintheorem}}
\newcommand{\emt}{\end{maintheorem}}
\newcommand{\bc}{\begin{corollary}}
\newcommand{\bl}{\begin{lemma}}
\newcommand{\ec}{\end{corollary}}
\newcommand{\el}{\end{lemma}}
\newcommand{\bo}{\begin{observation}}
\newcommand{\eo}{\end{observation}}
\newcommand{\bp}{\begin{proposition}}
\newcommand{\ep}{\end{proposition}}
\newcommand{\br}{\begin{remark}}
\newcommand{\er}{\end{remark}}
\newcommand{\bpr}{\begin{principle}}
\newcommand{\epr}{\end{principle}}

\def\PG{\mathbf{PG}}

\def\eop{\hspace*{\fill}$\blacksquare$}

\usepackage{amssymb} 
\usepackage{amsmath} 
\usepackage[utf8]{inputenc} 
\usepackage{latexsym}

\title {\bf On $k$-caps in $\PG(n, q)$, with $q$ even and $n \geq 4$}
\author {J. A. THAS \\Ghent University}
\address{Ghent University, Department of Mathematics, Krijgslaan 281, S25, B-9000 Ghent, Belgium}
\email{joseph.thas@ugent.be}
\begin {document}
\maketitle
\begin{abstract}
Let $m_2(n, q)$, $n \geq 3$, be the maximum size of $k$ for which there exists a complete $k$-cap in $\PG(n, q)$. In this paper the known bounds for $m_2(n, q), n \geq 4$, $q$ even and $q \geq 2048$, will be considerably improved. \\

\end{abstract}

\section{Introduction}
A {\em k-arc } of $\PG(2, q)$ is a set of $k$ points, no three of which are collinear; a {\em k-cap} of $\PG(n, q)$, $n \ge 3$, is a set of $k$ points, no three of which are collinear. A $k$-arc or $k$-cap is {\em complete} if it is not contained in a ($k+1$)-arc or ($k+1$)-cap. The largest value of $k$ for which a $k$-arc of $\PG(2, q)$, or a $k$-cap of $\PG(n, q)$ with $n\ge 3$, exists is denoted by $m_2(n, q)$. The size of the second largest complete $k$-arc of $\PG(2, q)$ or $k$-cap of $\PG(n, q)$, $n\ge 3$, is denoted by $m^\prime_2(n, q)$. \\

For any $k$-arc $K$ in $\PG(2, q)$ or $k$-cap $K$ in $\PG(n, q)$, $n \geq 3$, a {\em tangent} of $K$ is a line which has exactly one point in common with $K$. Let $t$ be the number of tangents of $K$ through a point $P$ of $K$ and let $\sigma_1(Q)$ be the number of tangents of $K$ through a point $Q \notin K$. Then for a $k$-arc $K$ $t + k = q + 2$ and for a $k$-cap K $t + k = q^{n-1} + q^{n-2} + \cdots + q + 2 $. \\

\begin{theorem}[\cite {JWPH: 87}]
\label{thm 1.1} 
If $K$ is a complete $k$-arc in $\PG(2, q)$, $q$ even, or a complete $k$-cap in $\PG(n, q)$, $n \geq 3$ and $q$ even, then $\sigma_1(Q) \leq t$ for each point $Q$ not in $K$.
\end{theorem}

\begin{theorem}
\label{thm1.2}
\begin {itemize}
\item [{\rm (i)}] $m_2(2, q) = q+2$, $q$ even  \cite {JWPH: 98};
\item [(ii)] $m_2(3, q) = q^2+1$, $q$ even, $q > 2$ \cite {JWPH: 85, RCB: 47, BQ: 52};
\item [(iii)] $m_2(n, 2) = 2^n$ \cite {RCB: 47};
\item [(iv)] $m_2(4, 4) = 41$ \cite {YE: 99};
\item [(v)] $m^\prime_2(n, 2) = 2^{n-1} + 2^{n-3}$ \cite {AAD: 90};
\item [(vi)] $m^\prime_2(3, 4) = 14$ \cite {JWPH: 87}.
\end {itemize}
\end{theorem}

\begin{theorem} [\cite {BS: 87, JAT: 87, JWPH: 98}]
\label{thm 1.3} 
Let $K$ be a $k$-arc of $\PG(2, q)$, $q$ even and $q > 2$, with $q - \sqrt{q} + 1 < k \le {q+1}$. Then $K$ can be uniquely extended to a ($q+2$)-arc of $\PG(2, q)$. 
\end{theorem}

The following result is the Main Theorem of \cite {JAT: 20...}.

\begin{theorem} [\cite {JAT: 20...}]
\label{thm 1.4}
\begin{equation}
m^\prime_2(3, q) < q^2 - (\sqrt{5} - 1)q + 5, q \ \text{even}\ , q \geq 8.
\end{equation}
\end{theorem}

As a corollary new bounds for $m_2(n, q)$, $q$ even, $q \ge 8$ and $n \ge 4$, are obtained.

\begin{theorem} [\cite {JAT: 20...}]
\label{thm 1.5}
\begin {itemize}
\item [{(i)}] $m_2(4, 8) \ge 479$;
\item [{(ii)}] for $q$ even, $q > 8$, \\ 
$m_2(4, q) < q^3 - q^2 + 2\sqrt{5}q - 8$;
\item [{(iii)}] $m_2(n, 8) \le 478.8^{n - 4} - 2(8^{n-5} + \cdots + 8 + 1) + 1, n \ge 5$;
\item [{(iv)}] for $q$ even, $q > 8, n \ge 5$, \\
$m_2(n, q) < q^{n -1} - q^{n - 2} + 2\sqrt{5}q^{n - 3} - 9q^{n - 4} - 2(q^{n - 5} + \cdots + q + 1) + 1$.
\end{itemize}
\end{theorem}

Combining the main theorem of \cite {SS: 93} with Theorem 1.4, there is an immediate improvement of the upper bound for $m^\prime_2(3, q), q \ge 2048$. This important remark is due to T. Sz\H{o}nyi.

\begin{theorem} [\cite {JAT: 20...}]
\label{thm 1.6} 
\begin{equation}
m^\prime_2(3, q) < q^2 - 2q + 3\sqrt{q} + 2, q \ \text{even}\ , q \ge 2048.
\end{equation}
\end{theorem}

Relying on Theorem 1.6, in the underlying paper new bounds for $m_2(n, q)$, $q$ even, $q \ge 2048, n \ge 4$, will be obtained.

\begin{theorem}
\label{thm 1.7}
For $q$ even, $q \ge 2048$,
\begin {itemize}
\item [{(i)}] $m_2(4, q) < q^3 - 2q^2 + 3q\sqrt{q} + 8q - 9\sqrt{q} - 6$,
\item [{(ii)}] $m_2(n, q) < q^{n - 1} - 2q^{n - 2} + 3q^{n - 3}\sqrt{q} + 8q^{n - 3} - 9q^{n - 4}\sqrt{q} - 7q^{n - 4} \\ - 2(q^{n - 5} + \cdots + q + 1) + 1, n \ge 5$.
\end{itemize}
\end{theorem}

\section{New bound for $m_2(4, q)$}

\begin{theorem}
\label{thm 2.1}
For $q$ even, $q \ge 2048$, 
\begin{equation}
m_2(4, q) < q^3 - 2q^2 + 3q\sqrt{q} + 8q - 9\sqrt{q} - 6.
\end{equation}
\end{theorem}
{\em Proof}. \quad
Assume, by way of contradiction, that $K$ is a complete $k$-cap of $\PG(4, q)$ with 
\begin{equation}
k \ge q^3 - 2q^2 + 3q\sqrt{q} + 8q - 9\sqrt{q} - 6.
\end{equation}
Then
\begin{equation}
t \le q^3 + q^2 + q + 2 - q^3 + 2q^2 - 3q\sqrt{q} - 8q + 9\sqrt{q} + 6,
\end{equation}
so
\begin{equation}
t \le 3q^2 - 3q\sqrt{q} - 7q + 9\sqrt{q} +8.
\end{equation}
We obtain a contradiction in several stages.\\

(\rm I) {\bf $K$ contains no plane $q$-arc}

Assume, by way of contradiction, that $\pi$ is a plane with $\vert\pi\cap K \vert = q$; let $\pi\cap K = Q$.

(\rm a) {\bf Suppose that $\delta_1, \delta_2, ..., \delta_5$ are distinct hyperplanes containing $\pi$, such that
\begin{equation}
\vert \delta_i \cap K \vert \ge q^2 - 2q + 3\sqrt{q} + 2, i = 1, 2, ..., 5
\end{equation}}

By Theorem 1.6 each $\delta_i \cap K$ can be extended to an ovoid $O_i$ of $\delta_i, i = 1, 2, ..., 5$. Hence $O_i \cap \pi$ is a $(q + 1)$-arc $Q \cup \lbrace N_i \rbrace, i = 1, 2, ..., 5$. Since $Q$ is contained in two $(q + 1)$-arcs at least three of the points $N_i$ coincide, say $N_1 = N_2 = N_3$. The joins of $N_1$ to the points of $\delta_i \cap K$, with $i = 1, 2, 3$, are tangents of $K$. \\
Hence
\begin{equation}
\sigma_1(N_!) \ge 3(q^2 - 3q + 3\sqrt{q} + 2) + q,
\end{equation}
so
\begin{equation}
\sigma_1(N_1) \ge 3q^2 - 8q + 9\sqrt{q} + 6.
\end{equation}
As $K$ is complete, $\sigma_1(N_1) \le t$. So 
\begin{equation}
3q^2 - 8q +9\sqrt{q} + 6 \le 3q^2 - 3q\sqrt{q} - 7q + 9\sqrt{q} + 8,
\end{equation}
that is
\begin{equation}
3q\sqrt{q} - q - 2 \le 0,
\end{equation}
clearly a contradiction.

(\rm b) {\bf Assume that there are at most 4 hyperplanes $\delta$ of $\PG(4, q)$ containing $\pi$ with $\vert \delta \cap K \vert \ge q^2 - 2q + 3\sqrt{q} + 2$} 

Then counting points of $K$ in hyperplanes containing $\pi$ gives
\begin{equation}
k < (q - 3)(q^2 - 3q + 3\sqrt{q} + 2) + 4(q^2 - q) + q,
\end{equation}
that is,
\begin{equation}
k < q^3 - 2q^2 + 3q\sqrt{q} + 8q - 9\sqrt{q} - 6.
\end{equation}
(Remark that any hyperplane containing $\pi$, has at most $q^2$ points in common with $K$.)

But $k \ge q^3 - 2q^2 + 3q\sqrt{q} + 8q - 9\sqrt{q} - 6$, clearly a contradiction.\\

(\rm II) {\bf There exists no hyperplane $\delta$ of $\PG(4, q)$ such that 
\begin{equation}
q^2 + 1 > \vert \delta \cap K \vert \ge q^2 - 2q + 3\sqrt{q} + 2
\end{equation}}
Suppose, by way of contradiction, that such a $\delta$ exists. Let $\delta \cap K = K^\prime$. Then $K^\prime$ can be extended to an ovoid $O$ of $\delta$. Let $N \in O \setminus K^\prime$ and let $N^\prime \in K^\prime$. Consider the $q + 1$ planes of $\delta$ containing the line $NN^\prime$. Each of these planes meets $O$ in a $(q + 1)$-arc, so by I each such plane meets $K^\prime$ in at most a $(q -1)$-arc.

Assume, by way of contradiction, that none of these intersections is a $(q - 1)$-arc. Counting the points of $K^\prime$ on these $q + 1$ planes gives 
\begin{equation}
\vert K^\prime \vert \le (q + 1)(q - 3) + 1,
\end{equation}
so
\begin{equation}
\vert K^\prime \vert \le q^2 - 2q - 2.
\end{equation}
As $\vert K^\prime \vert \ge q^2 - 2q + 3\sqrt{q} + 2$, there arises $3\sqrt{q} + 4 \le 0$, a contradiction.

So we may assume that $\vert \pi \cap K^\prime \vert = q - 1, \pi  \subset \delta, NN^\prime \subset \pi$. Consider all hyperplanes of $\PG(4, q)$ containing the plane $\pi$. Let $\theta$ be the number of such hyperplanes $\pi^\prime$ for which 
\begin{equation}
\vert \pi^\prime \cap K \vert \ge q^2 - 2q + 3\sqrt{q} + 2.
\end{equation}
By assumption $\theta \ge 1$.

First assume $\theta \ge 4$, hence there are at least 4 hyperplanes $\pi^\prime_1, \pi^\prime_2, \pi^\prime_3, \pi^\prime_4$ containing $\pi$ such that
\begin{equation}
\vert \pi^\prime_i \cap K \vert \ge q^2 - 2q +3\sqrt{q} +2.
\end{equation}
Consequently $\pi^\prime_i \cap K$ can be extended to an ovoid $O_i$ of $\pi^\prime_i$, with $i = 1, 2, 3, 4$. It follows that $O_i \cap \pi$ is a $(q + 1)$-arc $(\pi \cap K) \cup \lbrace N^\prime_i, N^{\prime\prime}_i \rbrace, i = 1, 2, 3, 4$. The $(q - 1)$-arc $\pi \cap K$ is extendable to a unique $(q + 2)$-arc $R$ of $\pi$, and each $(q + 1)$-arc of $\pi$ containing $\pi \cap K$ belongs to $R$ \cite {JWPH: 98} . So $\pi \cap K$ is contained in exactly 3 $(q + 1)$-arcs of $\pi$. It follows that there is at least one point $N$ which belongs to 3 of the 4 pairs $\lbrace N^\prime_i, N^{\prime\prime}_i \rbrace$. So the number of tangents $\sigma_1(N)$ of $K$ containing $N$ is at least
\begin{equation}
3(q^2 - 2q + 3\sqrt{q} + 2 - q + 1) + q - 1 = 3q^2 - 8q + 9\sqrt{q} + 8.
\end{equation}
As $\sigma_1(N) \le t$, there arises
\begin{equation}
3q^2 - 3q\sqrt{q} - 7q + 9\sqrt{q} + 8 \ge t \ge \sigma_1(N) \ge 3q^2 - 8q + 9\sqrt{q} + 8,
\end{equation}
so
\begin{equation}
3q\sqrt{q} - q \le 0,
\end{equation}
a contradiction.

Finally, assume $\theta \le 3$. Counting the points of $K$ in the $q + 1$ hyperplanes containing $\pi$, we obtain 
\begin{equation}
k < (q - 2)(q^2 - 2q + 3\sqrt{q} + 2 - q + 1) + 3(q^2 - q + 1) + q - 1,
\end{equation}
so
\begin{equation}
k < q^3 - 2q^2 + 3q\sqrt{q} + 7q - 6\sqrt{q} - 4.
\end{equation}
As
\begin{equation}
k \ge q^3 - 2q^2 + 3q\sqrt{q} + 8q - 9\sqrt{q} - 6,
\end{equation}
there arises 
\begin{equation}
q - 3\sqrt{q} - 2 < 0,
\end{equation}
a final contradiction. \\

(\rm III) {\bf For a point $N$ not on $K$, there do not exist planes $\pi_1$ and $\pi_2$ such that $\pi_1 \cap \pi_2 = \lbrace N \rbrace$ and such that $\pi_i \cap K$ is a $(q + 1)$-arc with nucleus $N, i = 1, 2$} 

Suppose, by way of contradiction, that such planes $\pi_1, \pi_2$ exist. Let $\delta$ be a hyperplane containing $\pi_1$. Then $\delta \cap K$ contains the $q + 1$ tangents of $\pi_1 \cap K$ through $N$ and one tangent of $\pi_2 \cap K$ through $N$. So $\delta \cap K$ has at least $q + 2$ tangents through $N$. Hence $\vert \delta \cap K \vert < q^2 + 1$.

Suppose that 
\begin{equation}
\vert \delta \cap K \vert < q^2 - 2q + 3\sqrt{q} + 2
\end{equation}
for any such hyperplane $\delta$. Counting points of $K$ in hyperplanes containing $\pi_1$ gives
\begin{equation}
k < (q + 1)(q^2 - 3q + 3\sqrt{q} + 1) + q + 1,
\end{equation}
so
\begin{equation}
k < q^3 - 2q^2 + 3q\sqrt{q} - q + 3\sqrt{q} + 2.
\end{equation}
As $k \ge q^3 - 2q^2 + 3q\sqrt{q} + 8q - 9\sqrt{q} - 6$, there arises a contradiction.

Consequently there exists a hyperplane $\delta$ through $\pi_1$ for which 
\begin{equation}
q^2 + 1 > \vert \delta \cap K \vert \ge q^2 - 2q + 3\sqrt{q} + 2,
\end{equation}
contradicting II. \\

(\rm IV) {\bf The tangents of $K$ through any point $N$ not on $K$ lie in a hyperplane}

Let $\delta$ be a hyperplane not containing $N$ and let $\mathcal V$ be the set of the intersections of $\delta$ with all tangents of $K$ through $N$. We will show that each point of $\mathcal V$ is on at least two lines contained in $\mathcal V$.

Let $R \in \mathcal V$ and let $r = RN$. Assume, by way of contradiction, that for at most one plane $\pi$ containing $r$ we have $\vert \pi \cap  K \vert \ge q - 1$. So
\begin{equation}
k \le (q^2 + q)(q - 3) + (q + 1),
\end{equation}
that is, 
\begin{equation}
k \le q^3 - 2q^2 - 2q + 1.
\end{equation}
 As $k \ge q^3 - 2q^2 + 3q\sqrt{q} + 8q - 9\sqrt{q} - 6$, there arises $3q\sqrt{q} + 10q - 9\sqrt{q} - 7 \le 0$, a contradiction. 
 
 Hence we may assume that for distinct planes $\pi, \pi^\prime$ containing $r$ we have
 \begin{equation}
 \vert \pi \cap K \vert, \vert \pi^\prime \cap K \vert \in \lbrace q - 1, q + 1 \rbrace.
 \end{equation}
 (By I no plane intersects $K$ in a $q$-arc.)
 
 We distinguish two cases.
 
 (\rm a) {\bf At least one of the planes $\pi, \pi^\prime$ intersects $K$ in a $(q - 1)$-arc}
 
 Say $\vert \pi \cap K \vert = q - 1$. Assume, by way of contradiction, that for no hyperplane $\delta^\prime$ containing $\pi$ we have $\vert \delta^\prime \cap K \vert = q^2 + 1$. Counting points of $K$ in hyperplanes containing $\pi$ gives by II
 \begin{equation}
 k < (q + 1)(q^2 - 2q + 3\sqrt{q} + 2 - q + 1) + q - 1,
 \end{equation}
 that is,
 \begin{equation}
 k < q^3 - 2q^2 + 3q\sqrt{q} + q + 3\sqrt{q} + 2.
 \end{equation}
 As $k \ge q^3 - 2q^2 + 3q\sqrt{q} + 8q - 9\sqrt{q} - 6$, there arises $7q - 12\sqrt{q} - 8 < 0$, a contradiction.
 
 So for at least one hyperplane $\delta^\prime$ containing $\pi$ we have $\vert \delta^\prime \cap K \vert = q^2 + 1$. But then $\vert \pi \cap K \vert = q + 1$, again a contradiction.
 
 (\rm b) {\bf $\vert \pi \cap K \vert = \vert \pi^\prime \cap K \vert = q + 1$}
 
 If $N$ is the nucleus of both $\pi \cap K$ and $\pi^\prime \cap K$, then there are two lines of $\mathcal V$ through $R$, namely $\pi \cap \delta$ and $\pi^\prime \cap \delta$.
 
 Therefore suppose that $N$ is not the nucleus of $\pi \cap K$. If for at most one hyperplane $\delta^\prime$ containing $\pi$ we have $\vert \delta^\prime \cap K \vert = q^2 + 1$, then counting points of $K$ in hyperplanes containing $\pi$ gives 
 \begin{equation}
 k < q^2 + 1 + q(q^2 - 3q + 3\sqrt{q} + 1),
 \end{equation}
 so
 \begin{equation}
 k < q^3 - 2q^2 + 3q\sqrt{q} + q + 1.
 \end{equation}
 As $k \ge q^3 - 2q^2 + 3q\sqrt{q} + 8q - 9\sqrt{q} - 6$, there arises $7q - 9\sqrt{q} - 7 < 0$, a contradiction.
 
 Consequently there are at least two hyperplanes $\delta_1$ and $\delta_2$ containing $\pi$ for which $\delta_i \cap K = O_i$ is an ovoid, $i = 1, 2$. then there is a plane $\pi_i$ of $\delta_i$ containing $N$ such that $N$ is the nucleus of the $(q + 1)$-arc $\pi_i \cap O_i = K_i, i = 1, 2.$ As $N$ is not the nucleus of $\pi \cap K$, we have $\pi \ne \pi_1 \ne \pi_2 \ne \pi$. The tangents of $K_i$ (which contain $N$) meet $\delta$ in the points of a line $l_i$ containing $R$, with $i = 1, 2$, and $l_1 \ne l_2.$.
 
 Consequently each point of $\mathcal V$ is on at least two lines contained in $\mathcal V$.
 
 If there existed two skew lines in $\mathcal V$, there would be two planes $\pi^\prime_1$ and $\pi^\prime_2$ on $N$, with $\pi^\prime_1 \cap \pi^\prime_2 = \lbrace N \rbrace$ and $N$ the nucleus of the $(q + 1)$-arcs $\pi^\prime_1 \cap K$ and $\pi^\prime_2 \cap K$. This is in contradiction with III. It follows that the lines of $\mathcal V$ all have a common point or all lie in a common plane. As each point of $\mathcal V$ is on at least two lines of $\mathcal V$, all lines of $\mathcal V$ lie in a plane. Hence $\mathcal V$ is subset of a plane, and so all tangents of $K$ containing $N$ lie in a hyperplane. \\
 
 (\rm V) {\bf The final contradiction}
 
 The final contradiction will be obtained by counting all tangents of $K$.
 
 Consider the function
 \begin{equation}
 G(x) = x(q^3 + q^2 + q + 2 - x).
 \end{equation}
 It attains its maximum value for 
 \begin{equation}
 x = \frac{1}{2} (q^3 + q^2 + q + 2).
 \end{equation}
 We have 
 \begin{equation}
 q^3 > k \ge q^3 - 2q^2 + 3q\sqrt{q} + 8q - 9\sqrt{q} - 6 > \frac{1}{2} (q^3 + q^2 + q + 2),
 \end{equation}
 so
 \begin{equation}
 kt = k(q^3 + q^2 + q + 2 - k) = G(k) > G(q^3) = q^3(q^2 + q + 2).
 \end{equation}
 
 All tangents containing a point $N$ not on $K$ lie in a hyperplane, which contains at most $q^2 + 1$ points of $K$. An ovoid of a hyperplane containing $N$ has exactly $q + 1$ tangents containing $N$. Hence $N$ is contained in at most $q^2$ tangents of $K$.
 
 Counting the pairs $(N, l)$, with $N \notin K$, $l$ a tangent of $K$ containing $N$, there arises
 \begin{equation}
 (q^4 + q^3 + q^2 + q + 1 - k)q^2 \ge ktq,
 \end{equation}
 so
 \begin{equation}
 (q^4 + q^3 + q^2 + q + 1 - q^3 + 2q^2 - 3q\sqrt{q} - 8q + 9\sqrt{q} + 6)q \ge kt,
 \end{equation}
 so
 \begin{equation}
 (q^4 + 3q^2 - 3q\sqrt{q} - 7q + 9\sqrt{q} + 7)q \ge kt > q^3(q^2 + q + 2),
 \end{equation}
 that is, 
 \begin{equation}
 q^3 - q^2 + 3q\sqrt{q} + 7q - 9\sqrt{q} - 7 < 0,
 \end{equation}
 a contradiction. \eop \\

 \section{New bound for $m_2(n, q), n \ge 5$}
 
 \begin{theorem}
 \label{thm 3.1}
 For $q$ even, $q \ge 2048, n \ge 5$
 \begin{equation}
 m_2(n, q) < q^{n - 1} -2q^{n - 2} + 3q^{n - 3}\sqrt{q} + 8q^{n - 3} - 9q^{n - 4}\sqrt{q} - 7q^{n - 4} - 2(q^{n - 5} +  \cdots + q + 1) + 1.
 \end{equation}
 \end{theorem}
 {\em Proof} \quad
 By 6.14(ii) of \cite{JWPH: 16} for $n \ge 5$ and $q > 2$, we have
 \begin{equation}
 m_2(n, q) \le q^{n - 4}m_2(4, q) - q^{n - 4} - 2(q^{n - 5} + \cdots + q + 1) + 1.
 \end{equation}
 From Theorem 2.1 the result follows. \eop \\


\begin{thebibliography}{99}
 \bibitem{RCB: 47}
 R. C. Bose, Mathematical theory of the symmetrical factorial design, $Sanky\bar{a}$ 8 (1947), 107 - 166.
 \bibitem{AAD: 90}
 A. A. Davydov and L. M. Tombak, Quasiperfect linear binary codes with distance 4 and complete caps in projective geometry, {\em Problems Inform. Transmission} 25 (1990), 265 - 275.
 \bibitem{YE: 99}
 Y. Edel and J. Bierbrauer, 41 is the largest size of a cap in $\PG(4, 4)$, {\em Des. Codes Cryptogr.} 16 (1999), 151 - 160.
 \bibitem{JWPH: 85}
 J. W. P. Hirschfeld, {\em Finite Projective Spaces of Three Dimensions}, Oxford University Press, Oxford, 1985, x + 316 pp.
 \bibitem{JWPH: 98}
 J. W. P. Hirschfeld, {\em Projective Geometries over Finite Fields}, Second Edition, Oxford University Press, Oxford, 1998, xiv + 555 pp.
 \bibitem{JWPH: 87}
 J. W. P. Hirschfeld and J. A. Thas, Linear independence in finite spaces, {\em Geom. Dedicata} 23 (1987), 15 - 31.
 \bibitem{JWPH: 16}
 J. W. P. Hirschfeld and J. A. Thas, {\em General Galois Geometries}, Second Edition, Springer, London, 2016, xvi + 409 pp.
 \bibitem{BQ: 52}
 B. Qvist, Some remarks concerning curves of the second degree in a finite plane, {\em Ann. Acad. Sci. Fenn. Ser. A} 134 (1952), 27 pp.
 \bibitem{BS: 87}
 B. Segre, Introduction to Galois geometries, {\em Atti Accad. Naz. Lincei Mem.} 8 (1987), 133 - 236 (edited by J. W. P. Hirschfeld).
 \bibitem{SS: 93}
 L. Storme and T. Sz\H{o}nyi, Caps in $\PG(n, q)$, $q$ even, $n \geq 3$, {\em Geom. Dedicata} 45 (1993), 163 - 169.
 \bibitem{JAT: 87}
 J. A. Thas, Complete arcs and algebraic curves in $\PG(2, q)$, {\em J. Algebra} 106 (1987), 451 - 464.
 \bibitem{JAT: 20...}
 J. A. Thas, On $k$-caps in $\PG(n, q)$, with $q$ even and $n \geq 3$, {\em Discrete Math.}, to appear.
 
 \end{thebibliography}
 \end {document}